\documentclass[12pt]{article}
\usepackage{amsmath}
\usepackage{amssymb}
\usepackage{amsthm}
\usepackage{amscd}
\usepackage{amsfonts}
\usepackage{graphicx}%
\usepackage{fancyhdr}

\theoremstyle{plain} \numberwithin{equation}{section}
\newtheorem{theorem}{Theorem}[section]
\newtheorem{corollary}[theorem]{Corollary}

\newtheorem{lemma}[theorem]{Lemma}

\theoremstyle{definition}
\newtheorem{definition}[theorem]{Definition}

\newtheorem{remark}[theorem]{Remark}

\def\Scd{\mathcal{S}(\D)}
\def\Scb{\mathcal{S}(\DT)}

\def\Scn{\mathcal{S}(\D^n)}

\def\INd{\mathcal{I}_N(\D)}
\def\INb{\mathcal{I}_N(\DT)}

\def\INn{\mathcal{I}_N(\D^n)}

\newcommand{\In}[1] { \mathcal{I}_N(\D^{#1})  }

\def\={\ = \ }

\def\l{\lambda}

\def\li{\lambda_i}

\def\={ = }

\def\C{\mathbb C}

\def\T{\mathbb T}
\def\D{\mathbb D}
\def\DT{\mathbb D^2} 

\def\be{\setcounter{equation}{\value{theorem}} \begin{equation}}
\def\ee{\end{equation} \addtocounter{theorem}{1}}
\def\beq{\begin{eqnarray*}}
\def\eeq{\end{eqnarray*}}

\def\bp{{\sc Proof: }}
\def\ep{{}{\hfill $\Box$} \vskip 5pt \par}

\def\bl{\begin{lemma}}
\def\el{\end{lemma}}
\def\bt{\begin{theorem}}
\def\et{\end{theorem}}
\def\bprop{\begin{prop}}
\def\eprop{\end{prop}}
\def\bd{\begin{definition}}
\def\ed{\end{definition}}
\def\br{\begin{remark}}
\def\er{\end{remark}}
\def\bexer{\begin{exercise}}
\def\eexer{\end{exercise}}

\newcommand{\Frac}[2] {\displaystyle{\frac{#1}{#2}} }

\newcommand{\oo}{\omega}

\title{A uniqueness theorem for bounded analytic functions on the polydisc}
\author{David Scheinker}
\date{}                                           

\begin{document}

\maketitle

\begin{abstract}
For each $n,N\geq1$ let $\INn$ denote the set of rational inner 
functions on $\D^n$ of degree strictly less than $N$.
We construct a set of points $\l_1,...,\l_{N^n}\in\D^n$ 
with the following property: if $f\in\INn$ and an analytic function $g$ 
maps $\D^n$ to $\D$ 
and satisfies $g(\li)=f(\li)$ for each $i=1,...,N^n$, then $g=f$ on $\D^n$.
In terms of the Pick problem on $\D^n$, our result implies 
that if $f\in\INn$, then the Pick problem with data $\l_1,...,\l_{N^n}$ and 
$f(\l_1),...,f(\l_{N^n})$ has a unique solution.
\end{abstract}

\section{Introduction}

Let $\D$ denote the unit disc in $\C$ and let $\T=\partial \D$. Let 
$\D^n$ and $\T^n$ denote the cartesian products of n copies of $\D$ 
and $\T$ in $\C^n$, respectively. A rational function $f$ on $\D^n$ is 
called \textbf{inner} if $f$ is analytic and $|f(\tau)|=1$ 
for almost every $\tau\in\T^n$. For a rational $f$, let 
$f=\frac{q}{r}$ for $q,r\in\C[z_1,...,z_n]$ relatively prime and 
define the \textbf{degree} of $f$ as the degree of $q$. 

Throughout this paper $N$ will denote a positive integer. Let $\INn$ denote the set of rational inner functions on $\D^n$ 
of degree strictly less than $N$. Let $\Scn$ denote the \textbf{Schur class} 
of $\D^n$, the set of analytic functions mapping $\D^n$ to $\D$. Our main result is that for each $N$ there exist points $\l_1,...,\l_{N^n}\in\D^n$ such that each $f\in \INn$ is uniquely determined in $\Scn$ by its values on $\l_1,...,\l_{N^n}$. In terms of the Pick problem on $\D^n$, our result implies 
that if $f\in \INn$, then the Pick problem with data $\l_1,...,\l_{N^n}$ and 
$f(\l_1),...,f(\l_{N^n})$ has a unique solution.

The following is our main result.

\pagebreak
\bt
\label{PLn}
For each $N$ there exist points $\l_1,...,\l_{N^{n}}\in\D^n$ with 
the following property: for each $f\in \INn$, if $g\in\Scn$ satisfies $g(\li)=f(\li)$ 
for $i=1,...,N^n$, then $g=f$ on $\D^n$.\\
Furthermore, the points $\l_1,...,\l_{N^n}$ may be chosen as follows.\\
Let $M=N^{n-1}$, for $r=2,...,n$ let $\tau_{1}^r,...,\tau_{{N}}^r\in\T$ 
be distinct and let $D_1,...,D_M$ be distinct analytic discs given by 
\[D_k:\D\to\D^n \textrm{ with } D_k(z)=(z,\tau_{i_{2,k}}^2 z...,\tau_{i_{n,k}}^n z).\]
If for each $k$, the points $\l_{k_1},...,\l_{k_N}\in D_k(\D)$ are distinct, 
then the points $\{\l_{k_j}\}$ for $k=1,...,M$ and $j=1,...,N$ have the above property. 
\et

This work began as an investigation of the connections between the 
Pick problem, rational inner functions and algebraic varieties on $\DT$, 
first discovered by Agler and McCarthy in \cite{agmc_dv}. I'm thankful to 
Jim Agler for his generous help in improving the exposition of this paper.

This paper is organized as follows. In section 2 we give background results on the Pick problem and prove the case $n=1$ of Theorem \ref{PLn}. In section 3 we give background results about rational inner functions on $\D^n$. In section 4 we prove a result on $\DT$ that we will use in the proof of Theorem \ref{PLn}. In section 5 we prove Theorem \ref{PLn}. In section 6 we prove a refinement of Theorem \ref{PLn} that yields stronger results for singular inner functions.

\section{The Pick problem on $\D^n$}

The \textbf{Pick problem} on $\D^n$ is to determine, given N distinct nodes\\
$\l_1,...,\l_N\in \D^n$ and $N$ target points $\oo_1,...,\oo_N \in \D$, whether there is an analytic function $f\in\Scn$ that satisfies $f(\li)=\oo_i$ for each $i=1,...,N$. If such an $f$ exists we call $f$ a \textbf{solution}. For $n=1$, Pick proved the following.

\bt{(Pick, 1916 \cite{pi16})}
\label{Pick2}
Fix a Pick problem on $\D$ with data $\l_1,..,\l_M$ and $\oo_1,...,\oo_N$. The following are equivalent.\\
\textbf{a.} The problem has a solution and the solution is unique.\\
\textbf{b.} The matrix  $P =  \displaystyle{ \left( \frac{ 1-\omega_i \overline{\omega_j}  } {  1-\l_i \overline{\l_j}   } \right)_{i,j=1}^N}$ is positive semi-definite and singular.\\
\textbf{c.} The problem has a rational inner solution $f$ and $\deg(f)<N$.
\end{theorem}

Part \textbf{b} of Theorem \ref{Pick2} will not be used in the present 
work but is included to make the theorem look more recognizable to 
those familiar with the result. The following lemma  is equivalent to the implication \textbf{c} $\to$ \textbf{a} in Theorem \ref{Pick2}. The lemma is case $n=1$ of 
Theorem \ref{PLn}.

\begin{lemma}{(Theorem \ref{PLn} case $n=1$)}\\
Fix $N>1$, $\l_1,...\l_N\in \D$ and $f\in\INd$.\\
If $g\in\Scd$ satisfies $g(\li)=f(\li)$ for $i=1,...,N$, then $g=f$ on $\D$.
\label{PL}
\end{lemma}

\section{Rational Inner functions on $\D^n$}

In \cite{rud69}, Rudin proved the 
following result about the structure of rational inner functions on $\D^n$.

\bt
\label{Rudin}
(Rudin 1969) Every $f\in\INn$ can be written
\[\varphi(z_1,...,z_n)=z_1^{d_1}\cdot\cdot\cdot z_n^{d_n}\Frac{\overline{q(\frac{1}{\overline{z_1}},...,\frac{1}{\overline{z_n}})}}{q(z_1,...,z_n)}\] for some polynomial $q(z_1,...,z_n)$ that does not vanish on $\D^n$.
\et
\noindent
We introduce two definitions before we employ Rudin's result.
\bd
For $n>m\geq 1$, we call an analytic function $E:\D^m\to\D^n$ an 
\textbf{analytic m-disc}. We use $E(\D^m)$ 
to denote the range of $E$. 
\ed
\noindent
\bd
Let $f\in\Scn$, $\tau\in\T$ and $E$ an analytic (n-1)-disc
\be
\textrm{given by }E:\D^{n-1}\to\D^n \textrm{ with } E(z_1,...,z_{n-1})=(z_1,...,z_{n-1},\tau z_1)
\label{DiscPar}
\ee
We define $f_E$ as follows
\[  f_E(z_1,...,z_{n-1})=f(E(z_1,...,z_{n-1})).\]
\ed
\noindent
The function $f_E$ is in $\mathcal{S}(\D^{n-1})$ and parametrizes 
the restriction of $f$ to $E$.

\begin{corollary}
\label{N1DiscDegree}
For $n\geq m>1$ and $\tau\in\T$, if $f\in \INn$ and\\
$E(z_1,...,z_{n-1})=(z_1,...,z_{n-1},\tau z_1)$, then $f_E \in \In{n-1}$.
\end{corollary}
\noindent
\bp
Since $f$ is inner, the denominator 
of $f$ has a non-zero constant term and Theorem \ref{Rudin} implies 
that $f$ can be written as follows.
\be f(z_1,...,z_n)=\tau \Frac{z^{d_1}\cdot\cdot\cdot z^{d_n}+r_0(z_1,...,z_n)}{1+q_0(z_1,...,z_n)} \textrm{ for some } \tau\in\T
\label{InnerForm}
\ee
where the degree of $f$ equals $d_1+...+d_n$ and each term of $r_0$ has degree less than or equal to $d_i$ in each $z_i$ and less than $d_i$ in at least one $z_i$. The corollary follows from substituting $f(z_1,...,z_{n-1},\tau z_1)$ into equation \ref{InnerForm}.
\ep

\section{A result on $\DT$}

We will use  the case $n=2$ of Lemma \ref{N1Disc} in the proof of Theorem \ref{PLn}.

We will use the following technical lemma to prove the case $n=2$ of Lemma \ref{N1Disc}. We use $B_\epsilon(z)$ to denote the ball of radius $\epsilon$ around $z$ and we use $m_{t,a}(z)$ to denote the automorphism of $\D$ given by $t\frac{z-a}{1-\bar a z}$.

\bl
Let $\tau_1,...,\tau_N\in\T$ be distinct and $E_1,...,E_M$ be analytic discs
\[ \textrm{ given by }E_i:\D\to\DT \textrm{ with } E_i(z)=(z,\tau_i z).\]
There exist $\tau\in\T$ and $\epsilon>0$ such 
that for every $t\in B_\epsilon(\tau)\cap\T$ and $a\in B_\epsilon(0)$,\\
the image of the analytic disc $C_{m_{t,a}}$ given by
\[ C_{m_{t,a}}:\D\to\DT \textrm{ with } C_{m_{t,a}}(z)=(z, m_{t,a}(z))\]
intersects each $E_i(\D)$ at a distinct point $(r_i,\tau_i r_i)$.\\
Furthermore, $C$, defined as the union of every $C_{m_{t,a}}(\D)$ over $t\in B_\epsilon(\tau)\cap\T$ and $a\in B_\epsilon(0)$
is a set of uniqueness for analytic functions on $\DT$.
\label{Mobius}
\el

\noindent
\bp
Fix $\tau\in\T$ such that $\tau \neq\tau_i$ for each $i$ and let 
$\epsilon_1>0$ be small enough so that for each $i$, 
$\tau_i\not \in B_{\epsilon_1}(\tau)$.
Let $C_m=C_{m_{\tau,a}}$, with $a$ to be specified later. 
The sets $C_m(\D)$ and $E_i(\D)$ intersect if and only if 
one of the roots of the equation $\tau_i z=m_{t,a}(z)$ lies in $\D$.
Let $r_i$ and $s_i$ denote the roots. If $a=0$ then $r_i=0$ and 
$s_i=\infty$ for each $i$. For sufficiently small $\epsilon_1>\epsilon>0$, 
if $a$ is perturbed 
away from zero and remains in $B_\epsilon(0)$, then each of the 
roots $r_i$ becomes non-zero and stays in $\D$.
That the roots $r_1,...,r_M$ are distinct follows from that they are non-zero and that $\tau_i \neq \tau_j$.

To see that $C$ is a set of uniqueness let $f$ be analytic on $\DT$ and 
suppose that $f|_C=0$. Fix $x\in\D$, $a\in B_\epsilon(0)$ and let 
\[A_x=\{(x,m_{t,a}(x))\in\DT: t\in B_\epsilon(\tau)\cap\T\} \subset C.\] 
Since $f(x,z)$ is an analytic function in 
the single variable $z$ and vanishes on the arc $A_x$, $f=0$. 
Since $f(x,\cdot)=0$ for each $x\in\D$, $f=0$ on $\DT$.
\ep

If Theorem \ref{PLn} holds for $n$ then the following lemma immediately 
follows for $n$. We prove the following lemma for $n=2$. 

\begin{lemma}
\label{N1Disc}
Fix $N$, let $f\in \INn$, let $\tau_1,...,\tau_N\in\T$ be distinct and let $E_1,...,E_N$ be analytic (n-1)-discs given by
\[
E_k:\D^{n-1}\to\D^n \textrm{ with } E_k(z_1,...,z_{n-1})=(z_1,...,z_{n-1},\tau_k z_1)
\]
If $g\in\Scn$ satisfies $g=f$ on each $E_k(\D^{n-1})$, then $g=f$ on $\D^n$.
\end{lemma}

{\sc Proof of lemma~\ref{N1Disc}}(case n=2):
By Lemma \ref{Mobius} there exists an analytic disc $C_m(\D)$
that intersects each of $E_1(\D),...,E_N(\D)$ at a distinct point\\
$R_i=(r_i,\tau_i r_i)$. Fix $f\in\INb$ and assume that $g\in\Scb$ satisfies $g=f$ on each 
$E_k(\D^{n-1})$. Let $f_m= f_{C_m}$ and $g_m= g_{C_m}$. Notice 
that $g_m\in\Scd$ and by Lemma \ref{N1DiscDegree}, $f_m\in \INd$. It 
follows that for  $i=1,...,N$,
\[g_m(r_i)=g(D_i(r_i))=g(r_i,\tau_i r_i)=f(r_i,\tau_i r_i)=f(D_i(r_i))=f_m(r_i)\]
Since $g_m(r_i)=f_m(r_i)$ for $i=1,...,N$, lemma \ref{PL} implies that 
$g_m=f_m$ on $\D$ and thus, $g=f$ on each $C_m(\D)$. By 
Lemma \ref{Mobius}, the discs $C_m(\D)$ sweep out a set of uniqueness and thus, $g=f$ on $\D^2$.
\ep

\section{Proof of Theorem \ref{PLn}}

In this section we use induction to prove Theorem \ref{PLn}. The 
case $n=1$ is Lemma \ref{PL}. Fix $n\geq 2$ and suppose that 
Theorem \ref{PLn} holds for each $m<n$. We show that Theorem 
\ref{PLn} holds holds for $n$ in 3 steps.

In the first step we fix $N$, fix a set of analytic (n-1)-discs\\ 
$E_1,...,E_{N}$, and fix a set of $N^{n-1}$ points 
$\{\l_{j_s}\}\subset\D^{n-1}$ to which we will imply the induction 
hypothesis. We lift the set $\{\l_{j_s}\}$ to the set of $N^n$ points 
$\{\l_{k_{j_s}}\}$ in $\D^n$ by letting $\l_{k_{j_s}}=E_k(\l_{j_s})$.

In the second step we apply the induction hypothesis to
show that for each $f\in\INn$, if $g\in\Scn$ satisfies 
$g(\l_{k_{j_s}})=f(\l_{k_{j_s}})$ for $k,j,s$, then $g=f$ on $E_1,...,E_{N}$.

In the third step we use Lemma \ref{N1Disc} (which holds for $n-1$ by the induction hypothesis) to show that since $g$ equals $f$ on $E_1,...,E_{N}$,
 $g=f$ on $\D^n$.\\

\noindent
\noindent
\textbf{STEP 1:} Fix $N$ and let $\tau_1,...,\tau_N\in\T$ be distinct and $E_1,...,E_{N}$ be analytic (n-1)-discs given by
\[E_k:\D^{n-1}\to\D^n \textrm{ with } E_k=(z_1,...,z_{m},\tau_i z_1).\]
Let $M=N^{n-2}$. For each $r=2,...,n-1$ let $\tau_{1}^r,...,\tau_{{N}}^r\in\T$ be distinct. Let $D_1,...,D_{N^{n-2}}$ be the $N^{n-2}$ analytic discs given by
\[ D_j:\D\to\D^{n-1} \textrm{ with } D_j(z)=(z,\tau_{i_{2,j}}^2 z...,\tau_{i_{n-1,j}}^{n-1} z).\]
For each $j$, let $\l_{j_1},...,\l_{j_N}\in D_j(\D) \subset \D^{n-1}$ be distinct and
lift each point $\l_{j_s}$ to $\D^n$, $N$ times, by letting $\l_{k_{j_s}}=E_k(\l_{j_s})$.\\

\noindent
\textbf{STEP 2:} Fix $f\in\INn$ and suppose $g\in\Scn$ satisfies $g(\l_{k_{j_s}})=f(\l_{k_{j_s}})$ for 
each $k,j,s$. For each $k$, let $f_k= f_{E_k}$ and $g_k= g_{E_K}$. Notice that\\ 
$g_k\in\mathcal{S}(\D^{n-1})$ and by Lemma \ref{N1DiscDegree}, $f_k\in\In{n-1}$. It follows that for\\
$k=1,...N$, $j=1,...,N^{n-2}$ and $s=1,...,N$,
\[g_k(\l_{j,s})=g(E_k(\l_{j,s}))=g(\l_{k,j,s})=f(\l_{k,j,s})=f(E_k(\l_{l,s}))=f_k(\l_{j,s}).\]
Since for each $k$, $g_k(\l_{j_s})=f_k(\l_{j_s})$ for each $j$ and $s$, the induction hypothesis implies that $g_k=f_k$ on $\D^{n-1}$. Thus, $g=f$ on each $E_k$.\\

\noindent
\textbf{STEP 3:} 
\noindent
If $n=2$, then case $n=2$ of Lemma \ref{N1Disc} implies that 
$g=f$ on $\DT$. Suppose $n\geq 3$.

For $\rho\in \T$ let $C_\rho$ be the analytic (n-1)-disc given by
\[C_\rho:\D^{n-1}\to\D^n \textrm{ with } C_\rho(z_1,...,z_{n-1})=(z_1,...,z_{n-2},z_{n-1},\bar \rho z_{n-1}).\]
\noindent
For each $\rho$, let $f_\rho= f_{C_\rho}$, $g_\rho= g_{C_\rho}$. Let  
$I_{\rho,k}:\D^{n-2}\to\D^n$ and \\
$H_{\rho,k}:\D^{n-2}\to\D^{n-1}$ be analytic (n-2)-discs such that 
\[I_{\rho,k}(\D^{n-2})=C_\rho(\D^{n-1})\cap E_k(\D^{n-1}) \textrm{ and }
H_{\rho,k}(\D^{n-2})=C_\rho^{-1}(I_{\rho,k}(\D^{n-2})).\]   
Since $g=f$ on $I_{\rho,1}(\D^{n-2}),...,I_{\rho,N}(\D^{n-2})$ 
it follows that $g_\rho=f_\rho$ on\\
$H_{\rho,1}(\D^{n-2}),...,H_{\rho,N}(\D^{n-2})$ and Lemma \ref{N1Disc} (which holds for $n-1$ by the induction hypothesis) implies that $g_\rho=f_\rho$. Thus, $g=f$ on $C_\rho$ and since 
$\D^n=\bigcup_{\rho\in\T}C_\rho$, it follows that $g=f$ on $\D^n$.
\ep

\section{A refinement and a question}

We call a rational inner function $f$ on $\D^n$ \textbf{singular}, if $f$ has 
a singular point on $\T^n$. In this section we show how Theorem \ref{PLn} may be refined to yield stronger results for singular inner functions.

For an analytic disc $D$ of the form in Theorem \ref{PLn}, let $\deg_D(f)$ equal the number of zeros of $f$ on $D$. Plugging $f(D(z))$ into formula \ref{InnerForm} implies that $\deg_D(f)$ is less than or equal to the degree of $f$. If $f$ has a singular point on $D$, 
then $\deg_D(f)$ is strictly less than the degree of $f$. Our proof of Theorem \ref{PLn} actually established 
the following refined theorem.
\bt If in Theorem \ref{PLn} the condition that there lie $N$
points on each analytic disc $D_k$ is replaced with the condition 
that there lie \\
$N_k=\deg_{D_k}(f)+1$ points on each $D_k$, then the conclusion still holds.
\label{PLnRefined}
\et


\begin{thebibliography}{10}


\bibitem{agmc_dv}
J.~Agler and J.E. M\raise.45ex\hbox{c}Carthy.
\newblock Distinguished Varieties.
\newblock {\em Acta Math.}, 194:133--153, 2005.

\bibitem{pi16}
G. Pick,
\newblock {\"U}ber die {Beschr\"ankungen} analytischer {F}unktionen, welche durch vorgegebene {F}unktionswerte bewirkt werden.
\newblock {\em Math. Ann.}, 77:7--23, 1916.

\bibitem{rud69}
W.~Rudin.
\newblock {\em Function Theory in {Polydiscs}}.
\newblock Benjamin, New York, 1969.

\bibitem{ds10}
D.~Scheinker.
\newblock Rational inner functions on the bidisc.
\newblock To appear.



\end{thebibliography}
\end{document}